\documentclass[a4paper, 12pt,oneside,reqno]{amsart}
\usepackage[a4paper]{geometry}
\usepackage{graphicx}
\usepackage[matrix,arrow,curve,cmtip]{xy}
\usepackage{amsfonts}
\usepackage{hyperref}
\usepackage{txfonts}

\usepackage{amsmath,amssymb,amsthm}

\def\wt{\mathop {\fam 0 wt} \nolimits}
\def\oo#1{\mathbin {{}_{(#1)}}}
\def\oon#1{\mathbin {{\circ}_{(#1)}}}
\def\Nov{\mathrm{Nov}}
\def\NovConf{\mathrm{NovConf}}
\def\ComConf{\mathrm{ComConf}}

\numberwithin{equation}{section}

\newtheorem{theorem}{Theorem}[section]
\newtheorem{lemma}[theorem]{Lemma}
\newtheorem{proposition}[theorem]{Proposition}
\newtheorem{corollary}[theorem]{Corollary}

\theoremstyle{definition}
\newtheorem{definition}[theorem]{Definition}
\newtheorem{example}[theorem]{Example}

\title{Differential envelopes of Novikov conformal algebras}

\author{P. S. Kolesnikov$^{1)}$}
\author{A. A. Nesterenko$^{2)}$}

\address{$^{1)}$Sobolev Institute of Mathematics, Novosibirsk, Russia.}
\address{$^{2)}$Novosibirsk State University, Novosibirsk,  Russia.}


\begin{document}
\begin{abstract}
A Novikov conformal algebra
is a conformal algebra such that its coefficient algebra
is right-symmetric and left commutative (i.e., it is an ``ordinary''
Novikov algebra).
We prove that every Novikov conformal algebra with a uniformly bounded
locality function on a set of generators can be embedded into
a commutative conformal algebra with a derivation.
In particular, every finitely generated Novikov conformal algebra
has a commutative conformal differential envelope.
For infinitely generated algebras this statement is not true in general.
\end{abstract}

\maketitle

\section{Introduction}

The class of nonassociative algebras called Novikov algebras
was introduced in the route of study of Hamiltonian operators
in formal variational calculus by Gelfand and Dorfman \cite{GelDor}.
Later, the same identities in a similar framework appeared
in the paper by Balinski and Novikov \cite{BalNovHydro}
as a tool to describe generalized Poisson brackets in the
theory of partial differential equations of hydrodynamic type.
The name ``Novikov algebras'' was proposed by Osborn~\cite{OsbNov}.
 
\begin{definition}[\cite{BalNovHydro, GelDor}]\label{defn:Novikov}
A Novikov algebra is a linear space $V$~$\Bbbk $ equipped with
a bilinear operation $\circ $ satisfying the following
identities:
\begin{gather}
(x\circ y)\circ
z - x\circ (y\circ z) = (x\circ z)\circ y - x\circ (z\circ y),
                     \label{eq:NovikovID-RSym} \\
x\circ (y\circ z) = y\circ (x\circ z).
                     \label{eq:NovikovID-LCom}
\end{gather}
\end{definition}

The study of structure theory for Novikov algebras was initiated in \cite{Zel}:
Zelmanov described finite-dimensional simple Novikov algebras over a
field of characteristic zero. Further structure results were obtained,
in particular, by Osborn \cite{OsbSim, OsbInf} and Xu \cite{XuSimNov, XuClass}.
Among the recent studies on Novikov algebras we should mention
the results on solvability and nilpotency \cite{ShestZhang2020, ZhelUmir2021},
and on algebraic (in)dependence in Novikov algebras \cite{DuiUmir2021}.

One of the main constructions in the theory of Novikov algebras
is the following operation-changing functor (Gelfand construction)
from the category of commutative differential algebras to the
category of Novikov algebras.

\begin{example}\label{exmp:Derivation}
Let $A$ be an associative and commutative algebra
equipped with a derivation
$d:A\to A$.
Then the same space
$A$ equipped with the operation
 $a\circ b = d(a)b$, $a,b\in A$,
 is a Novikov algebra
denoted by~$A^{(d)}$.
\end{example}

In \cite{DL2002}, the free Novikov algebra generated by a set~$B$ was constructed.
Let us state those details of that construction that are essential for
subsequent exposition.
Denote by $B^{(\omega )}$ the disjoint union of countably many copies of~$B$:
\[
B^{(\omega )} = \{x^{(p)} \mid x\in B, p\in \mathbb Z_+\},
\]
where $\mathbb Z_+$ is the set of nonnegative integers.
Introduce the derivation $d$ on the polynomial algebra $\Bbbk [B^{(\omega )}]$
defined by its value on the generators:
$d(x^{(n)})=x^{(n+1)}$ for $x\in B$, $n\in \mathbb Z_+$.

\begin{definition}\label{defn:weight}
Define the {\em weight} $\wt(\cdot )$ of a monomial
in the variables $B^{(\omega )}$ as follows:
\[
\wt (x_1^{(k_1)} \dots x_l^{(k_l)}) = k_1+\dots + k_l-l.
\]
If all monomials that appear in a polynomial
$f\in \Bbbk [B^{(\omega )}]$
with nonzero coefficients have the same weight $w\in \mathbb Z$
then
$f$ is said to be  $\wt$-{\em homogeneous},
and
$\wt(f)=w$.
\end{definition}

It is easy to see that
$\wt(fg) = \wt(f)+\wt(g)$ and
$\wt(d(f))=\wt(f)+1$ for all
$\wt $-homogeneous polynomials
$f,g\in \Bbbk[B^{(\omega )}]$,
and
\[
\Bbbk[B^{(\omega )}]
= \bigoplus\limits_{w\in \mathbb Z} F_w,
\]
where $F_w$ consists of zero and all $\wt $-homogeneous polynomials of weight~$w$.
In particular, we have
$d(F_{-1})F_{-1}\subseteq
F_0F_{-1}\subseteq F_{-1}$.

\begin{example}\label{exmp:FreeNovikov}[\cite{DL2002}]
Let $V$ be the subalgebra of the Novikov algebra
$\Bbbk [B^{(\omega )}]^{(d)}$
generated by the set
$\{x^{(0)}\mid x\in B\}$ which is a copy of~$B$.
Then $V=F_{-1}$, and this is a free Novikov algebra generated by~$B$.
\end{example}

The first statement is relatively easy to prove.
For example, the monomial
$u = x^{(2)}y^{(0)}z^{(0)}$ of weight $-1$
can be presented as
\[
u = d(x^{(1)}y^{(0)})z^{(0)} - x^{(1)}y^{(1)}z^{(0)} =(x^{(0)}\circ
y^{(0)})\circ z^{(0)} - x^{(0)}\circ (y^{(0)}\circ z^{(0)}).
\]

The second statement in Example~\ref{exmp:FreeNovikov}
highly non-trivial. In the original proof \cite{DL2002}
the authors construct a linearly complete set of monomials in the
free Novikov algebra $\mathrm{Nov} \langle B\rangle $ which turns to be a
basis, an alternative proof (in the non-commutative setting) was
proposed in \cite{SartKol_EUR}.

It was shown in \cite{BokutChenZhang} (see also \cite{KolSart_UEXM})
that for every (not necessarily free) Novikov algebra
$V$ there exists a commutative algebra
$A$ with a derivation $d$ such that $V$
is isomorphic to a subalgebra of the Novikov algebra
$A^{(d)}$.

\begin{definition}[\cite{Kac1998}]\label{defn:Conformal}
Let $\Bbbk $ be a field of characteristic zero.
A {\em conformal algebra} is a linear space $C$ over $\Bbbk $,
equipped with a linear operator
$\partial:C\to C $ and bilinear operation
$(\cdot \oo\lambda \cdot)$ with a range in the space of polynomials
$\Bbbk [\partial,\lambda ]\otimes_{\Bbbk [\partial ]}C \simeq C[\lambda ]$
in a formal variable $\lambda $ such that
the following axioms hold:
\begin{equation}\label{eq:sesqui-lin}
(\partial
a\oo\lambda b ) = -\lambda (a\oo\lambda b),
\quad
\partial(a\oo\lambda  b) =
(a\oo\lambda \partial b) - \lambda(a \oo \lambda b),
\end{equation}
for all
$a,b\in C$.
\end{definition}

In other words, if $C$ is a conformal algebra then
for every $a,b\in C$
there is a uniquely defined polynomial
$(a\oo\lambda b)\in C[\lambda ]$.
The coefficients of $(a\oo\lambda b)$ at $\lambda ^n/n!$
are denoted by
$(a\oo n b)\in C$,
$n\in \mathbb Z_+$:
\begin{equation}\label{eq:lambda-oo}
(a\oo\lambda b) =
\sum\limits_{n\ge 0} \dfrac{\lambda^n}{n!}   (a\oo n b)\in C[\lambda ].
\end{equation}
Thus, conformal algebras may be considered as algebraic systems
with infinitely many ``ordinary'' bilinear operations
$(a,b)\mapsto
(a\oo n b)$, $n\in \mathbb Z_+$,
satisfying the {\em locality} propery:
for every $a,b\in C$ there exists
$N\in \mathbb Z_+$ such that $(a\oo n b) = 0$ for all
$n\ge N$.
The minimal such $N\ge 0$ is denoted $N_C(a,b)$,
the function $N_C$ is called {\em locality function} of a conformal algebra~$C$.

A derivation on a conformal algebra $C$ is a linear map
$D:C\to C$ such that
$D\partial = \partial D$ and
$D(a\oo\lambda b) = (D(a)\oo\lambda
b) + (a\oo\lambda D(b))$ for all $a,b\in C$.
For example, $D=\partial $ is always a derivation due to
\eqref{eq:sesqui-lin}.

\begin{definition}[\cite{Roitman1999, KacForDist}]
Given a conformal algebra $C$, one may construct its
{\em coefficient algebra} $\mathcal A(C)$ as follows.
For every $n\in \mathbb Z$ denote by $\hat{A}(n)$ an isomorphic
copy of the space $C$, an image of $a \in C$ in $\hat{A}(n)$ is denoted by $a(n)$.
Let $\hat{A} = \bigoplus_{n\in \mathbb Z} \hat A(n) $,
and let $E \subset \hat{A}$ be the subspace
generated by all elements of the form
\[
(\partial a)(n) + na(n-1),
\quad  a \in C, \ n \in \mathbb Z.
\]
Then the space $\mathcal A(C)$ is set to be the quotient
$\hat{A}/E$, and the product
\begin{equation}\label{eq:CoeffProduct} a(n)b(m) = \sum\limits_{s\ge
0} \binom{n}{s} (a\oo s b)(n+m-s).
\end{equation}
is a well-defined bilinear operation on $\mathcal A(C)$.
Hereinafter we identify $a(n)$ with $a(n)+E\in \mathcal A(C)$.
\end{definition}

The rule $C\mapsto \mathcal A(C)$ is a functor from the
category of conformal algebras to the category of ``ordinary'' algebras.
Moreover, every derivation $D$ on a conformal algebra $C$
induces a derivation $d$ on $\mathcal A(C)$
by the rule
\[
d(a(n)) = (Da)(n),\quad a\in C,\ n\in \mathbb Z.
\]

The importance of coefficient algebras in the theory of conformal algebras
is explained by the construction of a conformal algebra by means of
formal distributions (see \cite{KacForDist, Roitman1999}).
Let $C$ be a conformal algebra and let
$\mathcal A(C)$ be the coefficient algebra of $C$.
Consider the space of formal distributions
$\mathcal A(C)[[z,z^{-1}]]$, i.e., two-side infinite formal power
series in a variable $z$ with coefficients from $\mathcal A(C)$.
For every $a\in C$ construct the series
\[
a(z) = \sum\limits_{n\in \mathbb Z} a(n)
z^{-n-1} \in \mathcal A(C)[[z,z^{-1}]].
\]
The linear map $a\mapsto a(z)$ is injective,
and
 $\partial a \mapsto \dfrac{d}{dz}a(z)$.
Moreover,
\begin{equation}\label{eq:ResFormulaOPE}
(a\oo{n} b)(z) =
\mathop{\fam 0 Res}\limits_{w=0} a(w)b(z)(w-z)^n
\end{equation}
for all
$a,b\in C$,
$n\in \mathbb Z_+$.
Here $\mathop{\fam 0 Res}\limits_{w=0}$ stands for the residue
of a formal series at $w=0$, i.e., the coefficient at $w^{-1}$
(depending on~$z$).

The locality condition on $C$ is equivalent to
\[
a(w)b(z) (w-z)^{N_C(a,b)} = 0, \quad a,b\in C,
\]
in the space
$\mathcal A(C)[[z,z^{-1},w,w^{-1}]]$.
In the coefficient-wise  form, the latter equation is equivalent
to the following system of equations in $\mathcal A(C)$:
\begin{equation}\label{eq:LocalityCoeff}
\sum\limits_{s\ge 0} (-1)^s
\binom{N_C(a,b)}{s} a(n-s)b(m+s) = 0 ,
\end{equation}
for all $a,b\in C$,
$n,m\in \mathbb Z$.
Thus, if \eqref{eq:LocalityCoeff} holds for some
$a,b\in C$ and for all $n,m\in \mathbb Z$
then it remains valid when we replace $N_C(a,b)$ with an
integer~$N>N_C(a,b)$.

The functor $\mathcal A(\cdot )$ may also be applied
to split the entire class of conformal algebras into smaller
subclasses that play the role of varieties of algebras.

\begin{definition}[\cite{Roitman1999}]\label{defn:VarConformal}
Let $\mathfrak{U}$ be a class of ``ordinary'' algebras
 (e.g., associative, commutative, etc.).
A conformal algebra $C$ is said to be a $\mathfrak{U}$-conformal algebra
if the corresponding coefficient algebra $\mathcal A(C) $ belongs to~$\mathfrak{U}$.
\end{definition}

If $\mathfrak{U}$ is a variety defined by a family of identities
then the class of $\mathfrak{U}$-conformal algebras
may also be described in terms of identities that hold
for the operations $\partial $ and $(\cdot \oo\lambda \cdot )$.

\begin{example}[\cite{KacForDist}]
Let us state conformal versions of common identities.
In what follows, $\lambda $ and $\mu $ are independent variables,
and if we are given an expression $(a\oo\lambda b)$ in the form
\eqref{eq:lambda-oo} then $(a\oo{-\partial-\lambda} b)$
is obtained by replacing the variable
$\lambda $ with the operator
$-\partial-\lambda$.
\begin{itemize}
\item[1)]
Associativity:
\[
a\oo\lambda (b\oo\mu c) = (a\oo\lambda b)\oo{\lambda+\mu } c.
\]
\item[2)] (Anti-)commutativity:
\[
(a\oo\lambda b ) = \pm
(b\oo{-\partial-\lambda } a).
\]
\item[3)] Jacobi identity:
\[
a \oo \lambda (b \oo \mu c)
- b \oo \mu (a\oo \lambda c)
= (a \oo \lambda b) \oo {\lambda + \mu}c.
\]
\end{itemize}
\end{example}

In a similar way, the identities
\eqref{eq:NovikovID-RSym} and
\eqref{eq:NovikovID-LCom}
of Novikov algebras
may be transformed to their conformal versions.

\begin{definition}[\cite{HongLi2015,KolNest2023}] \label{defn:NovikovConformal}
A conformal algebra
$C$ with operations
$\partial $ and $(\cdot \oon \lambda \cdot )$
is a Novikov conformal algebra
if
\begin{gather}
(a\oo\lambda b)\oo{\lambda+\mu} c -
a\oo{\lambda } (b\oo{\mu } c) = (a\oo \lambda c)\oo{-\partial-\mu} b -
a\oo\lambda (c\oo{-\partial-\mu } b)  , \label{eq:RSymLambda}
\\
a\oo\lambda
(b\oo \mu c) = b\oo\mu (a\oo\lambda c), \label{eq:LComLambda}
\end{gather}
for all $a,b,c\in C$.
\end{definition}

\begin{example}\label{exmp:DerConf}
Let a commutative conformal algebra $C$ with a binary operation
$(\cdot\oo\lambda \cdot)$ be equipped with a derivation~$D$.
Then the space $C$ relative to the same operator $\partial $
and new binary operation
$(a\oon\lambda b) = (Da\oo\lambda b)$, $a,b\in C$,
is a Novikov conformal algebra denoted $C^{(D)}$.
\end{example}

Other examples of Novikov conformal algebras may be obtained
from (ordinary) Novikov--Poisson algebras.

\begin{definition}[\cite{Xu1997}]
A linear space $P$ equipped with two bilinear operations  $\circ$ and $*$,
is called a {\em Novikov--Poisson algebra} if the following conditions hold:
 $(P,\circ )$ is a Novikov algebra,
 $(P,*)$ is associative and commutative,
 and
\begin{equation}\label{eq:NP-additional}
(a\circ b)*c = a\circ (b*c),
\quad
(a*b)\circ c - a*(b\circ c) = (a*c)\circ b - a*(c\circ b)
\end{equation}
for all $a,b,c\in P$.
\end{definition}

\begin{example}[\cite{HongLi2015}]
Let $(P,\circ, *)$ be a Novikov--Poisson algebra.
Consider the free  $\Bbbk[\partial]$-module
$V = \Bbbk [\partial ]\otimes P$
equipped with an operation
$(\cdot \oon \lambda \cdot)$ defined by the rule
\[
(a\oon \lambda b ) = a \circ b + \lambda (a*b), \quad a,b\in P.
\]
(On the entire $\Bbbk [\partial ]$-module $V$ this operation extends uniquely
by means of \eqref{eq:sesqui-lin}.)
Then $V$ is a Novikov conformal algebra (such examples
are natural to call {\em quadratic} Novikov conformal algebras, they were
studied in
\cite{KolNest2023}).
\end{example}

In this paper we study the following natural question:
whether the construction of a Novikov conformal algebra from Example~\ref{exmp:DerConf}
is generic in the same sense as for ordinary Novikov algebras (Example~\ref{exmp:Derivation}).

\begin{definition}\label{defn:SpecialNovConf}
A Novikov conformal algebra
$V$ is said to be {\em special}
if there exists a commutative conformal algebra $C$
with a derivation $D$ such that $V$ is isomorphic to a subalgebra of $C^{(D)}$.
\end{definition}

In \cite{KolNest2023} it was shown that every
quadratic Novikov conformal algebra $V$
is special.
As a main result this paper, we prove that every Novikov conformal algebra
$V$ generated (as a conformal algebra) by a set $B$ such that
$N_V(a,b)$ is uniformly bounded for all $a,b\in B$
is special.
(In particular, every finitely generated Novikov conformal algebra is special.)
We also state an example to show that the embedding statement is not true
in general for conformal algebras: there exist non-special Novikov conformal algebras.

\section{Free Novikov conformal algebras}

The class of Novikov (as well as associative, commutative, Lie, etc.)
conformal algebras is not a variety in the sense of universal algebra:
an infinite Cartesian product of conformal algebras may not be a
conformal algebra due to locality issues.
However, one may construct ``relatively free'' conformal algebras
in certain classes by means of an additional restriction on the
locality of generators. For associative, commutative, and Lie
conformal algebras it was done in~\cite{Roitman1999}.
In a similar way, a family of free Novikov conformal algebras
may be constructed.

Let us fix a nonempty set $X$ and a function
$N: X\times X\to \mathbb Z_+$.
Consider the set
\[
B =X\times \mathbb Z = \{a(n) \mid a\in X, n\in \mathbb Z\} 
\]
and the free Novikov algebra
$\Nov\langle B\rangle $ generated by $B$ over a field $\Bbbk$
($\mathrm{char}\,\Bbbk =0$, as above).
The multiplication in $\Nov\langle B\rangle $
will be denoted by the symbol~$\circ $.

Denote by $\mathcal N = \mathcal N(X;N)$
the quotient algebra $\Nov\langle B\rangle/J(N) $,
where $J(N)$ is the ideal generated by all elements of the form
\[
\sum\limits_{s\ge 0} (-1)^s \binom{N(a,b)}{s}
 a(n-s)\circ b(m+s), \quad 
 a,b\in X,\ n,m\in \mathbb Z.
\]
Then formal distributions of the form
\begin{equation}\label{eq:N-generators}
a(z) = \sum\limits_{n\in \mathbb Z} a(n) z^{-n-1}\in 
\mathcal N[[z,z^{-1}]],\quad a\in X, 
\end{equation}
are pairwise local, i.e.,
\[
a(w)\circ b(z) (w-z)^{N(a,b)} = 0
\]
for all $a,b\in X$.

Recall the key statement which is essential for constructing
conformal algebras by means of formal distributions.

\begin{lemma}[Dong Lemma, c.f. \cite{HLie1997}]
Let $A$ be an associative or Lie algebra,
and let $a(z),b(z),c(z)\in A[[z,z^{-1}]]$ be formal distributions.
Suppose these three series are pairwise mutually local, i.e.,
for every $x,y\in \{a,b,c\}$ there exists $N(x,y)\in \mathbb Z_+$
such that
\[
 x(w)y(z)(w-z)^{N(x,y)}=0 \in A[[z,z^{-1},w,w^{-1}]].
\]
Then $(a\oo n b)(z)$ and $c(z)$ are also mutually local,
where the first series is given by \eqref{eq:ResFormulaOPE}.
\end{lemma}

Since the derivation $\partial = d/dz$ also preserves locality,
the Dong Lemma implies that every set of pairwise mutually local
formal distributions over an (associative or Lie) algebra $A$
generates a conformal algebra immersed into $A[[z,z^{-1}]]$.

It is shown in \cite{BokKol2024}
that the Dong Lemma remains valid for Novikov algebras.
Hence, the formal distributions \eqref{eq:N-generators}
generate a conformal algebra $\NovConf(X; N)$ in $\mathcal N[[z,z^{-1}]]$
relative to the operations
$\partial = d/dz$ and
$n$-products denoted
$(\cdot\oon n\cdot )$, $n\in \mathbb Z_+$.
The universal property of $\NovConf(X; N)$
is similar to whose in \cite{Roitman1999} for Lie or associative
conformal algebras.

Let $X$ and $N$ be as above.

\begin{theorem}
Let $V$ be a Novikov conformal algebra
generated by a set $X$ such that $N_V(a,b)\le N(a,b)$
for all $a,b\in X$.
Then there is a unique epimorphism of conformal algebras
\[
\psi : \NovConf(X;N) \to V
\]
such that
$\psi(a)=a$ for $a\in X$.
\end{theorem}

\begin{proof}
The coefficient algebra $\mathcal A(V)$
is a Novikov algebra generated by
$B=\{a(n)\in \mathcal A(V)\mid a\in X, n\in \mathbb Z \}$.
Due to the locality property the
following relations hold on $ \mathcal A(V)$:
\[
\sum\limits_{s\ge 0} (-1)^s \binom{N_V(a,b)}{s}
 a(n-s)\circ b(m+s) = 0, \quad 
 a,b\in X,\ n,m\in \mathbb Z.
\]
Since $N_V(a,b)\le N(a,b)$, we may replace $N_V$ with $N$
to get valid relations on
$\mathcal A(V)$.
Hence, there exists a homomorphism of Novikov algebras
$f: \mathcal N(X;N)\to \mathcal A(V)$,
such that $a(n) \in \mathcal N(X;N)$ maps to  $a(n) \in \mathcal A(V)$.
Let us extend $f$ in a coefficient-wise way to a linear map
\[
\psi : \mathcal N(X;N)[[z,z^{-1}]]\to \mathcal A(V)[[z,z^{-1}]].
\]
The map $\psi $ commutes with $d/dz$
and preserves conformal products defined by
\eqref{eq:ResFormulaOPE}
since the coefficients of such products may be expressed
via the coefficients of factors by means of addition and multiplication.
Note that
\[
\psi (a(z)) = a(z)\in \mathcal A(V)[[z,z^{-1}]], \quad a\in X.
\]
Hence, the image of
 $\NovConf(X;N)$ under $\psi $ belongs to
the conformal algebra generated by $a(z)$, $a\in X$,
inside the space $\mathcal A(V)[[z,z^{-1}]]$.
The latter conformal algebra is isomorphic to~$V$ \cite{Roitman1999}.

Therefore, the restriction of $\psi $ onto $\NovConf(X;N)$
is the desired homomorphism of conformal algebras.
Such a homomorphism is unique since the images of generators are fixed.
\end{proof}

The ``PBW-Theorem'' for Novikov algebras
proved in \cite{BokutChenZhang} allows us
to consider the Novikov algebra
$\mathcal N(X;N)$ as a subalgebra in the
differential commutative algebra constructed as follows.

Consider the set
\[
B^{(\omega )} = \{a^{(p)}(n) \mid a\in X,\ n\in \mathbb Z,\ p\in \mathbb Z_+\},
\]
we will identify the elements $a(n)\in B$ with
$a^{(0)}(n) \in B^{(\omega )}$.
Define a derivation $d$ on the (commutative) polynomial algebra
$\Bbbk [B^{(\omega )}]$ by means of its values on the variables:
\begin{equation}\label{eq:d-Derivation}
d(a^{(p)}(n)) = a^{(p+1)}(n).
\end{equation}
Denote by $I(N)$ the ideal of $\Bbbk [B^{(\omega )}]$,
generated by the polynomials
\begin{equation}\label{eq:Poly-1-generators}
f_{a,b}(n,m) = \sum\limits_{s\ge 0} (-1)^s\binom{N(a,b)}{s}
a^{(1)}(n-s)b(m+s), \quad a,b\in X,\ n,m\in \mathbb Z,
\end{equation}
and all their derivatives $d^s(f_{a,b}(n,m))$, $s\in \mathbb Z_+$.
Then
\[
U(X;N) = \Bbbk [B^{(\omega )}]/I(N)
\]
is a differential commutative algebra
with an induced derivation also denoted~$d$.
The image of $B$ is a subset of the Novikov algebra
$U(X; N)^{(d)}$, and the subalgebra generated by
all these $a(n)+I(N)$, $a(n)\in B$,
is isomorphic to $\mathcal N(X;N)$.
In other words,
\[
I(N) \cap F_{-1} = J(N),
\]
where $F_{-1}$ is the space of $\wt$-homogeneous polynomials in
the variables $B^{(\omega )}$ of weight $-1$,
and $J(N)$ in the right-hand side
is identified with its isomorphic image under the embedding
$\Nov\langle B\rangle \subset \Bbbk [B^{(\omega )}]$.
Since the image of the entire $\Nov\langle B\rangle$
under this embedding is equal to $F_{-1}$
\cite{DL2002},  we obtain
\begin{equation}\label{eq:NXB-weight}
\mathcal N(X;N) \simeq \{f+I(N) \mid 
f\in \Bbbk [B^{(\omega )}], \wt(f)=-1 \}\subset U(X;N)^{(d)}.
\end{equation}

\section{On the speciality of free Novikov conformal algebras}

As in the previous section,
let $X$ be a nonempty set,
$N: X\times X \to \mathbb Z_+$ be a fixed function,
the elements of
$B = X\times \mathbb Z$ are written in the form
$a(n)$, $a\in X$, $n\in \mathbb Z$.
Recall that
$B^{(\omega )}$ consists of the symbols
\[
a^{(p)}(n), \quad a\in X,\ n\in \mathbb Z, \ p\in \mathbb Z_+. 
\]
Denote by $F$ the (commutative) algebra of polynomials
$\Bbbk [B^{(\omega )}]$, 
and consider
$f_{a,b}(n,m) \in F$
given by \eqref{eq:Poly-1-generators}.
As above, let $I(N)$ be the differential ideal
of $F$ generated by all
$f_{a,b}(n,m)$ relative to the derivation \eqref{eq:d-Derivation}.

Similarly, denote
\[
X^{(\omega)} = \{a^{(p)} \mid a \in X, 
p \in \mathbb{Z}_{+} \}.
\]
Assume
$\bar N$  is a function $X^{(\omega)} \times X^{(\omega)}
\to \mathbb{Z}_{+}  $.
Denote
\begin{equation}\label{eq:f-pq}
f_{a,b}^{p,q}(n,m) =
\sum \limits_{t \ge 0} (-1) ^{t} 
 \binom{\bar N(a^{(p)}, b^{(q)})}{t} 
 a^{(p)}(n-t) b^{(q)}(m+t) \in F
\end{equation}
for $a,b\in X$, $p,q\in \mathbb Z_+$, $n,m\in \mathbb Z$.
Let $K(\bar N)$ stand for the ideal of $F$ generated by all polynomials
$f_{a,b}^{p,q}(n,m)$ and all their derivatives relative to the derivation
\eqref{eq:d-Derivation}.

Note that the polynomial
$f_{x,y}^{p,q}(n,m)$
is equal to the coefficient at
$z^{-m-1}w^{-n+\bar N(x^{(p)},y^{(q)})-1}$
in the formal power series
\[
x^{(p)}(w)y^{(q)}(z)(w-z)^{\bar N(x^{(p)}, y^{(q)})} 
 \in F[[z,z^{-1},w,w^{-1}]].
\]
Here, as above,
$x^{(p)}(z) =\sum\limits_{s\in \mathbb Z} x^{(p)}(s)z^{-s-1}\in F[[z,z^{-1}]]$. 

The main technical problem in this section is resolved by
the following statement.

\begin{proposition}\label{prop:LocalityExists}
Suppose there exists a constant $M\in \mathbb Z_+$
such that $N(a,b)\le M$ for all $a,b\in X$.
Then one may find a funstion
$\bar N: X^{(\omega)} \times X^{(\omega)}
\to \mathbb{Z}_{+}  $
such that
\begin{equation}\label{eq:IdealCaps}
F_{-1} \cap I(N) = F_{-1} \cap K(\bar N).
\end{equation}
\end{proposition}

\begin{proof}
Let us set
\[
\bar N(a^{(p)}, b^{(q)}) = 
\begin{cases}
  3M, & p=q=0, \\
  N(a,b), & p=1, q=0, \\
  (p + q)M, & \text{otherwise}.
\end{cases}
\]
In particular, $\bar N(a^{(1)}, b^{(0)}) =  N(a,b)$, so
$f^{1,0}_{a,b}(n,m) = f_{a,b}(n,m)$. 
Therefore,
$ I(N)\subseteq K(\bar N)$
and
$F_{-1} \cap I(N) \subseteq F_{-1} \cap K(\bar N)$.
It remains to show the converse inclusion.

A generic polynomial $h$ in $K(\bar N)$ has the following form:
\begin{equation}\label{eq:Generic-J-element}
h= \sum\limits_i \alpha_i u_i d^{l_i} \big (
f_{x_i,y_i}^{p_i,q_i}(n_i,m_i)
\big ),
\end{equation}
where $\alpha_i\in \Bbbk $, $u_i$ are monomials from $F$,
$l_i\ge 0$.

Every polynomial
$d^{l_i}\big ( f_{x_i,y_i}^{p_i,q_i}(n_i,m_i) \big )$
is $\wt$-homogeneous of weight $l_i+p_i+q_i-2$,
so in the case $\wt(h)=-1$ for every summand in
\eqref{eq:Generic-J-element}
we have
\[
\wt(u_i) = 1-p_i-q_i-l_i.
\]
Let us show that every summand in
\eqref{eq:Generic-J-element} belongs to~$I(N)$.

Suppose
\[
h=u d^{l}  \big ( f_{x,y}^{p,q}(n,m) \big ), 
\quad \wt(u)=1-p-q-l.
\]
In order to prove $h\in I(N)$ let us
consider three cases:
\begin{enumerate}
    \item\label{Case:1} 
      $p=q=l=0$, then $\wt (u)=1$;
    \item\label{Case:2}
       $p+q+l=1$, then $\wt(u)=0$;
    \item\label{Case:3} 
        $p+q+l\ge 2$, then $\wt(u)<0$.
\end{enumerate}

{\sc Case~\ref{Case:1}.}
Since the weight of $u$ is positive, it has to contain a letter
$a^{(r)}(k)\in B^{(\omega )}$
for some $a\in X$, $k\in \mathbb Z$, $r\ge 2$.
Let us show
$a^{(r)}(k) f_{x,y}^{0,0}(n,m) \in I(N)$,
which would imply $h\in I(N)$.

First, note that
$a^{(r)}(k) f_{x,y}^{0,0}(n,m)$ is a coefficient of the
formal distribution
\[
a^{(r)}(\zeta ) x(w)y(z) (w-z)^{3M} \in F[[z,z^{-1}, w,w^{-1}, \zeta, \zeta^{-1}]].
\]
Proceed by induction on $r\ge 2$.
For $r=2$ (denote $a^{(2)}$ by $a''$; we will also
use the notation $d(f)=f'$ for a polynomial $f\in F$),
transform the latter distribution as follows:
\begin{multline}\label{eq:Series00}
a''(\zeta ) x(w)y(z) (w-z)^{3M}    
= a''(\zeta ) x(w)y(z) (w-\zeta +\zeta -z)^{2M} (w-z)^M \\
= a''(\zeta ) x(w)y(z) 
 \big ((w-\zeta)^M P +(\zeta -z)^{M} Q \big ) (w-z)^M \\
= a''(\zeta ) x(w)y(z) (w-\zeta)^M (w-z)^M P
  +a''(\zeta ) x(w)y(z)(\zeta -z)^{M} (w-z)^M Q,
\end{multline}
where $P,Q\in \Bbbk [z,w,\zeta]$ are those polynomials
that appear in the distribution of
$((w-\zeta) +(\zeta -z))^{2M}$ via the binomial formula, i.e.,
\[
(w-z)^{2M} = (w-\zeta)^M P +(\zeta -z)^{M} Q.
\]
Next, rewrite the summands in the right-hand side of
\eqref{eq:Series00}:
\begin{multline*}
a''(\zeta ) x(w)y(z) (w-\zeta)^M (w-z)^M   
=\big( 
  (a'(\zeta )x(w))' y(z) - a'(\zeta)x'(w)y(z) \big )
   (w-\zeta)^M (w-z)^M   \\
= \big (a'(\zeta )x(w)(w-\zeta)^M\big )' y(z) (w-z)^M
  - a'(\zeta)(w-\zeta)^M \big( x'(w)y(z) (w-z)^M\big)
\end{multline*}
Since $M\ge N(a,x), N(x,y)$,
every coefficient in both terms above contains a factor
which can be expressed
either via some derivatives of
$f_{a,x}(n_i,m_i)$ or via $f_{x,y}(n_i,m_i)$,
for some integer $n_i$, $m_i$.
Hence, all these coefficients belong to the ideal $I(N)$.
The multiplication by the polynomial $P$ with scalar coefficients
leads to linear combinations of the latter coefficients,
so they remain to be in $I(N)$.
Therefore, all coefficients of the first summand in the right-hand side of
\eqref{eq:Series00} belong to $I(N)$.

In a very similar way, all coefficients of the second summand
in \eqref{eq:Series00} are also in $I(N)$:
\begin{multline*}
a''(\zeta ) x(w)y(z) (\zeta-z)^M (w-z)^M   
=\big( 
  (a'(\zeta )y(z))' x(w) - a'(\zeta)y'(z)x(w) \big )
  (\zeta-z)^M (w-z)^M   \\
= \big (a'(\zeta )y(z)(\zeta-z)^M \big )' x(w) (w-z)^M
  - a'(\zeta)(\zeta-z)^M \big( y'(z)x(w) (w-z)^M\big).
\end{multline*}
Therefore, in particular,
 $a''(k) f_{x,y}^{0,0}(n,m) \in I(N)$
as desired.

Assume it is already proved that
$a^{(r-1)}(k) f_{x,y}^{0,0}(n,m)\in I(N)$
for some $r>2$ and for all
$k,n,m\in \mathbb Z$.
Then all coefficients of the distribution
$a^{(r-1)}(\zeta ) x(w)y(z) (w-z)^{3M}$
belong to $I(N)$ along with all their derivatives (relative to $d$).
Consider the following distribution with coefficients from $I(N)$:
\begin{multline*}
d\big( a^{(r-1)}(\zeta ) x(w)y(z) (w-z)^{3M} \big ) \\
= a^{(r)}(\zeta ) x(w)y(z) (w-z)^{3M}
+a^{(r-1)}(\zeta ) x'(w)y(z) (w-z)^{3M}
+a^{(r-1)}(\zeta ) x(w)y'(z) (w-z)^{3M}.
\end{multline*}
In the second and third summands of the right-hand side
all coefficients may be expressed via appropriate
$f_{x,y}(n_i,m_i)$ and $f_{y,x}(n_i,m_i)$
since $3M\ge N(x,y),N(y,x)$.
Hence, the coefficients of the first summand
also lie in $I(N)$ as desired.

{\sc Case \ref{Case:2}.}
Since $l,p,q\ge 0$, the equality $p+q+l=1$ implies either $l=0$, $p+q=1$
or $l=1$, $p=q=0$.

If $l=0$, $p=1$, $q=0$ then $f_{x,y}^{p,q}(n,m)$ coincides with
$f_{x,y}(n,m)$ from \eqref{eq:Poly-1-generators},
so $h\in I(N)$.

If $l=0$, $p=0$, $q=1$ then
$\bar N(x^{(0)}, y^{(1)}) = M$, and
replacing the index $t$ with $M-s$ in \eqref{eq:f-pq} we obtain
\[
f_{x,y}^{0,1}(n,m)
=\sum\limits_{s\ge 0} (-1)^{M-s} \binom{M}{M - s} y^{(1)}(m+M-s) x^{(0)} (n-M+s).
\]
This is an element from $I(N)$ since $M\ge N(y,x)$ by the construction.

Finally, if $d=1$, $p=q=0$ then
$d(f_{x,y}^{0,0} (n,m)) = f_{x,y}^{1,0} (n,m) + f_{x,y}^{0,1} (n,m)$,
where both summands are in $I(N)$.

Note that $I(N)$ is a differential ideal by definition, so
\[
d^l (f_{x,y}^{0,0}(n,m)) \in I(N)
\]
for all $l>0$.

{\sc Case \ref{Case:3}.}
Suppose
$h=u d^{l}  \big ( f_{x,y}^{p,q}(n,m) \big )
 \in K(\bar N)$, 
$\wt(u) = 1-p-q-l<0$, as above.

Proceed by induction on the pair
$(p+q,l)$ to show $h\in I(N)$.
The base of induction is given by the cases
$p=q=0$, $l\ge 0$
and
$p+q=1$, $l=0$, which are already considered.

Assume we have proved
\[
v d^r(f^{i,j}_{x,y}(n,m)) \in I(N)
\]
for every monomial $v\in F$,
$\wt(v) = 1-i-j-r$,
under the condition
$i+j<p+q$ or if $i+j=p+q$, $r<l$.
%

Case \ref{Case:3}(a): $l=0$.
Since $\wt(u)<0$, the word $u$ must contain a letter of negative weight,
so are only $a(k)$ for $a\in X$, $k\in \mathbb Z$.
Then $h$ is obtained from one of the coefficients
of the formal distribution
\[
a(\zeta )x^{(p)}(w)y^{(q)}(z) (w-z)^{\bar N(x^{(p)}, y^{(q)})}
\]
by a monomial $v\in F$ of weight $2-p-q$,
where $u=va(k)$.
Let us transform the latter distribution as follows:
\begin{equation}
    \label{eq:Series-pq}
a(\zeta )x^{(p)}(w)y^{(q)}(z) (w-z)^{(p+q)M}
=
a(\zeta )x^{(p)}(w)y^{(q)}(z) (w-\zeta)^{pM} P
+
a(\zeta )x^{(p)}(w)y^{(q)}(z) (\zeta-z)^{qM} Q,
\end{equation}
where the polynomials $P,Q\in \Bbbk [z,w,\zeta ]$ are obtained
from
\[
(w-z)^{(p+q)M}
= ( (w-\zeta) + (\zeta - z))^{pM+qM}
 = (w-\zeta)^{pM} P + (\zeta - z)^{qM} Q.
\]
Consider the summands in the right-hand side of \eqref{eq:Series-pq}.

If $p>0$ and $q>0$ then both these summands contain factors equal to
some coefficients of the distributions
$a(\zeta )x^{(p)}(w) (w-\zeta)^{pM}$
and
$a(\zeta ) y^{(q)}(z) (\zeta-z)^{qM}$,
respectively.
By the inductive assumption, the coefficients of
$a(\zeta )x^{(p)}(w) (w-\zeta)^{pM}$
drop to the ideal $I(N)$ when multiplied by
$p-1$ letters from $B$.
Similarly, the product of a coefficient of
$a(\zeta ) y^{(q)}(z) (\zeta-z)^{qM}$
and $q-1$ letters from $B$ also belongs to $I(N)$.
The monomial $v$ contains at least $p+q-1$ letters from $B$,
so $va(k)f_{x,y}^{p,q}(n,m)\in I(N)$.

If $p=0$ or $q=0$ then the expression \eqref{eq:Series-pq}
may be simplified.
For example, if $q=0$ then $p>1$, and
\begin{multline}\label{eq:Series-p0}
 a(\zeta )x^{(p)}(w)y(z) (w-z)^{pM}
 =  a(\zeta )\big( x^{(p-1)}(w)y(z) \big )' (w-z)^{pM}
 -  a(\zeta )x^{(p-1)}(w)y'(z) (w-z)^{pM} \\
 = a(\zeta )\big( x^{(p-1)}(w)y(z) (w-z)^{(p-1)M} 
  \big )'(w-z)^M 
 -  a(\zeta )\big( x^{(p-1)}(w)y'(z) (w-z)^{pM} \big).
\end{multline}
The coefficients of the first summand in the right-hand side of
\eqref{eq:Series-p0} are of the form
$a(k_i)d\big (f^{p-1,0}_{x,y}(n_i,m_i)\big )$.
By induction, assume all such polynomials drop into $I(N)$
when multiplied with a monomial $v$ of weight $2-p$.
The coefficients of the second summand are of the form
$a(k_i) f^{p-1,1}_{x,y}(n_i,m_i)$,
they also drop into $I(N)$ under multiplication by $v$,
since the case of positive indices $p$, $q$ is already considered.

The case $p=0$, $q>1$ is completely analogous.

Case \ref{Case:3}(b): $l>0$.
For brevity, denote
$f = f_{x,y}^{p,q}(n,m) \in F$.
We have to show
$h=ud^{l}(f) \in I(N)$ whenever
$\wt(u) = 1-p-q-l$. 

As above, present
$u = v a(k)$ for some
$a(k)\in B\subset B^{(\omega )}$.
Then $\wt(v) = \wt(u') = 2-p-q-l$,
and
\[
u d^l(f) = (u d^{l-1}(f))' - u'd^{l-1}(f)
=(a(k) v d^{l-1}(f) )' - u' d^{l-1}(f).
\]
By induction,
$v d^{l-1}(f)\in I(N)$ and $u' d^{l-1}(f)\in I(N)$,
so $h\in I(N)$ as desired.
\end{proof}

Let $\bar N$ be the function from Proposition \ref{prop:LocalityExists}.
Denote
\[
A(X;\bar N) = \Bbbk [B^{(\omega )}]/ K(\bar N),
\]
and let $d$ stand for the inherited derivation \eqref{eq:d-Derivation}
of the latter algebra.

\begin{corollary}\label{cor:Cor1}
The Novikov algebra $\mathcal N(X;N)$
is isomorphic to a subalgebra of
$A(X;\bar N)^{(d)}$ generated by~$B$.
\end{corollary}

\begin{proof}
As we have already mentioned
(see \eqref{eq:NXB-weight}), the Novikov algebra
$\mathcal N(X;B)$ is isomorphic to the image of the subspace
$F_{-1}$ in $U(X;N)^{(d)}$, where $U(X;N) = F/I(N)$.
In the proof of Proposition \ref{prop:LocalityExists}
we noted that $I(N)\subseteq K(\bar N)$.
Hence, there is a homomorphism of commutative algebras
$\varphi : U(X;N)\to A(X;\bar N)$ commuting with $d$, such that
\[
\varphi : f+I(N) \mapsto f+K(\bar N), \quad f\in F.
\]
If $f+I(N)\in \ker\varphi \cap \mathcal N(X;N)$
then $\wt(f)=-1$ and $f\in K(\bar N)$. By Proposition~\ref{prop:LocalityExists}
we conclude
$f\in F_{-1}\cap K(\bar N)\subseteq I(N)$, 
i.e., $\ker\varphi \cap \mathcal N(X;N) = 0$,
so the restriction of $\varphi $ on $\mathcal N(X;N)$ is injective.
\end{proof}

\begin{theorem}\label{thm:Embedding}
Let $X$ be a nonempty set,
and let $N: X\times X\to \mathbb Z_+$ be a uniformly bounded function,
i.e., there exists a constant
$M$ such that $N(a,b)\le M$ for all $a,b\in X$.
Then there exists a commutative conformal algebra
$C$
with a derivation $D$
such that $\NovConf(X;N)$ is isomorphic to a subalgebra of $C^{(D)}$.
\end{theorem}

\begin{proof}
Given $X$ and $N$, construct $\bar N$
as in Proposition~\ref{prop:LocalityExists},
and denote $A(X,\bar N)=A$ for brevity.
This is a commutative algebra with a derivation $d$.
Consider formal distributions of the form
\begin{equation}\label{eq:Comm-generators}
a^{(p)}(z) = \sum\limits_{n\in \mathbb Z} a^{(p)}(n) z^{-n-1} \in A[[z,z^{-1}]], 
\quad a\in X, \ p\in \mathbb Z_+.
\end{equation}
(Here we identify $a^{(p)}(n) \in B^{(\omega )}$
with their images in $A$.)
All series \eqref{eq:Comm-generators}
are pairwise local by the construction of the ideal $K(\bar N)$.
Hence, they generate a conformal algebra $C$ immersed into  $A[[z,z^{-1}]]$
relative to the operations \eqref{eq:ResFormulaOPE}.

Define a derivation $D: C \to C$ induced by $d$
in a coefficient-wise way:
\[
D: \sum\limits_{n\in \mathbb Z} f(n)z^{-n-1}
\mapsto \sum\limits_{n\in \mathbb Z} d(f(n))z^{-n-1},
\quad f(n)\in A.
\]
Indeed,
$D(a^{(p)}(z)) = a^{(p+1)}(z)$,
$D$ commutes with $\partial = d/dz$,
and for every $f(z),g(z)\in A[[z,z^{-1}]]$
we have
\begin{multline*}
D((f\oo n g)(z)) =
\mathop{\fam 0 Res}\limits_{w=0} D(f(w)g(z)(w-z)^n)
= \mathop{\fam 0 Res}\limits_{w=0} 
\big [ (Df(w)) g(z) (w-z)^n + f(w)D(g(z)) (w-z)^n\big ] \\
=(Df\oo n g)(z) + (f\oo n Dg)(z),
\end{multline*}
so $D(C)\subseteq C$.

Since
\[
(Df)(w) g(z) = f(w)\circ g(z) \in A^{(d)}[[w,w^{-1},z,z^{-1}]], 
\]
the Novikov conformal algebra
$C^{(D)}$ contains a subalgebra
generated by all distributions of the form
$a(z)=a^{(0)}(z)\in \mathcal N(X;N)[[z,z^{-1}]]$.
By Corollary~\ref{cor:Cor1}, $\mathcal N(X;N)[[z,z^{-1}]] \subset A[[z,z^{-1}]]$,
so by the construction this subalgebra is exactly
$\NovConf(X;N)$.
\end{proof}

\section{Speciality of homomorphic images and a non-special example}

Suppose $X$ is a nonempty set and $N: X\times X \to \mathbb Z_+$
is a uniformly bounded function, as above.

Note that the differential commutative conformal enveloping algebra $C$
of the free Novikov conformal algebra $\NovConf(X;N)$
constructed above is not
necessarily a free commutative
conformal algebra generated by
the set $X^{(\omega )}=\{a^{(p)}\mid a\in X, p\in \mathbb Z_+\}$
with the restriction $\bar N$ on the locality of generators.
Indeed, the construction of the algebra $A(X;\bar N)$
involves the derivatives of locality relations that do not appear
in the coefficient algebra of a free commutative conformal algebra.

In other words, let
$\ComConf(X^{(\omega )};\bar N)$ be the free commutative conformal algebra
generated by $X^{(\omega )}$ with respect to a locality
function $\bar N: X^{(\omega )}\times X^{(\omega )} \to \mathbb Z_+$.
Then the mapping $a^{(p)}\mapsto a^{(p+1)}$,
$a\in X$, $p\in \mathbb Z_+$,
may not be continued (in general) to a derivation of  $\ComConf(X^{(\omega )};\bar N)$
due to locality issues.

However, the commutative conformal algebra $C$ constructed in the proof
of Theorem~\ref{thm:Embedding} is a quotient of
$\ComConf(X^{(\omega )};\bar N)$ modulo an ideal $K$ generated by all derivatives
of the locality relations:
\begin{equation}\label{eq:generateK}
 \sum\limits_{s\ge 0} \binom{d}{s} a^{(p+s)}\oo n b^{(q+d-s)} = 0,
 \quad n\ge \bar N(a^{(p)}, b^{(q)}),\ d\ge 0,
\end{equation}
for all $a^{(p)}, b^{(q)}\in X^{(\omega )}$.

Define {\em weight} of a monomial in $\ComConf(X^{(\omega )};\bar N)$
in the same way as it was done in Definition~\ref{defn:weight}
for ordinary polynomials in $B^{(\omega }$:
\[
 \wt (a^{(p)})=p-1, \quad \wt (f\oo n g) = \wt(f)+\wt(g), \ \wt(\partial f)=\wt (f),
\]
for $a^{(p)}\in X^{(\omega )}$. This is a well-defined function
since the defining relations of $\ComConf(X^{(\omega )}, \bar N)$
considered as an abstract algebraic system with operations $\partial(\cdot)$
and $(\cdot \oo n \cdot)$, $n\in \mathbb Z_+$,
are $\wt$-homogeneous.

In particular, the elements \eqref{eq:generateK} are also $\wt $-homogeneous,
so we can define $\wt $ function on the conformal algebra $C=\ComConf(X^{(\omega )};\bar N)/K$.
By definition, $\wt(Du)=\wt(u)+1$ for a $\wt$-homogeneous element $u\in C$.

\begin{proposition}\label{prop:WeightNovConf}
 Suppose $f\in C$.
 Then $f\in \NovConf(X;N)$ if and only if $f$ is $\wt$-homogeneous and $\wt (f)=-1$.
\end{proposition}

The ``only if'' part is obvious from the definition of an operation in
$C^{(D)}$:
$(u\oon n v) = (Du\oo n v)$. If $u,v\in C$ are of weight $-1$ then
so are $(u\oon n v)$ and $\partial (u)$.
The ``if'' part may be proved in the very same way as it is done
for ordinary differential polynomials in \cite{DL2002}:
every formal expression in terms of $(\cdot \oo n \cdot )$ in the variables
$X^{(\omega )}$ may be written via the operations $(\cdot \oon n \cdot )$
in $X$.

\begin{theorem}
If $N$ is a uniformly bounded locality function on a non-empty set $X$
then every homomorphic image $V$ of the free Novikov conformal algebra
$\NovConf(X;N)$ is special.
\end{theorem}

\begin{proof}
 Let $I$ be an ideal of $\NovConf(X;N)$ such that $V = \NovConf(X;N)/I$.
Consider the set $X^{(\omega )}$ and the function $\bar N$ as in the proof
of Proposition~\ref{prop:LocalityExists}.
Then $\NovConf(X;N)$ embeds into the quotient
\[
 C = \ComConf(X^{(\omega )};\bar N)/K,
\]
where $K$ is the differential ideal generated by \eqref{eq:generateK}.
In terms of Proposition~\ref{prop:WeightNovConf}
we can say that $K$ contains no nonzero elements of weight $-1$.

Denote by $(I)_D$ the ideal of $C$ generated by all derivatives of $I$.
Then $(I)_D \cap \NovConf(X;N) = I$.
Indeed, an element from $(I)_D$ is a sum of expressions like
\[
h= \partial ^s (x_1 \oo{n_1} x_2 \oo{n_2} \dots x_k\oo{n_k} D^p(f)),
 \quad x_i\in X^{(\omega)}, \ f\in I,
\]
with some bracketing, $n_i, s, p\in \mathbb Z_+$.
If $\wt(h)=-1$ then $h\in I$ since it can be expressed via $x_i$ and $f$
$(\cdot \oon n \cdot)$ and $\partial $.

As a result, $\NovConf(X;N)/I$ is isomorphic to a subalgebra of
$(C/(I)_D)^D$, as desired.
\end{proof}

On a finite number of generators, a locality function is uniformly bounded.
Hence, we have

\begin{corollary}
A finitely generated Novikov conformal algebra is special.
\end{corollary}

In order to show that the condition $N(X,X)\le M$ is essential,
let us state an example of a non-special Novikov conformal algebra.

Let
\[
X = \{x, v_0, v_1, v_2, \ldots \}
\]
be a countable alphabet, and let
$W$
be the free
$\Bbbk [\partial ]$-module generated by
$X$ equipped with the following conformal product:
\begin{equation}\label{eq:Operation-Exmp}
(v_k\oon\lambda x) = (\partial+\lambda )^k v_k,
\quad
(x\oon \lambda x) = (x\oon\lambda v_k) = (v_k\oon\lambda v_p)=0
\end{equation}
for $k,p\in \mathbb Z_+$.

In order to show that $W$ is a Novikov conformal algebra
it is enough to check the relations
\eqref{eq:RSymLambda} and \eqref{eq:LComLambda}
for $a,b,c\in X$.
Note that the operation $(\cdot \oon\lambda \cdot )$
defined by  \eqref{eq:Operation-Exmp}
has the following property:
all expressions of the form
$a\oon\lambda (b\oon \mu c)$
are zero. Hence, \eqref{eq:LComLambda} holds trivially.
The remaining terms of \eqref{eq:RSymLambda} are:
\[
(a\oon\lambda b)\oon{\lambda+\mu} c
=
(a\oon \lambda c)\oon{-\partial-\mu} b .
\]
The only case when these expressions are nonzero is
$a=v_k$, $b=c=x$. Then
\[
(a\oon\lambda b)\oon{\lambda+\mu} c
=( (\partial +\lambda)^k v_k \oon{\lambda +\mu} x)
= (-\mu)^k (\partial+\lambda+\mu)^k v_k
\]
by \eqref{eq:sesqui-lin}.
On the other hand,
\begin{multline*}
(a\oon \lambda c)\oon{-\partial-\mu} b
=
( (\partial +\lambda)^k v_k \oon{-\partial-\mu} x)
=(\partial+\mu+\lambda)^k (v_k \oon{-\partial-\mu} x) \\
= (\partial+\mu+\lambda)^k (\partial + (-\partial-\mu))^k
v_k =
(\partial+\mu+\lambda)^k (-\mu)^k v_k.
\end{multline*}
Hence, \eqref{eq:RSymLambda} also holds on $W$.

\begin{proposition}
The Novikov conformal algebra $W$ cannot be embedded into
 $C^{(D)}$
for neither commutative conformal algebra
$C$ with a derivation~$D$.
\end{proposition}

\begin{proof}
Assume the converse: there exists a commutative conformal algebra
$C$  with conformal multiplication
$(\cdot\oo\lambda \cdot )$ and with a derivation $D:C\to C$
such that
$W\subseteq C$,
$(u\oon \lambda v) = (Du \oo \lambda v)$
for $u,v \in W$.

Denote by $N_x$ the value of the locality function of $x\in W$ with itself in $C$:
$N_x = N_C(x,x) \in \mathbb Z_+$.
Choose an integer $n\ge N_x$ and compute
\[
W\ni ((v_k\oon 0 x)\oon n x)
=
(D(D(v_k)\oo 0 x)\oo n x)
=((D^2v_k \oo 0 x)\oo n x) + ((Dv_k\oo 0 Dx)\oo n x)
\]
for every $k\ge 0$.
Both summands in the right-hand side of the last expression are zero:
\[
((D^2v_k \oo 0 x)\oo n x) = (D^2v_k \oo 0 (x\oo n x))=0,
\]
\[
((Dv_k\oo 0 Dx)\oo n x)=(Dv_k\oo 0 (Dx \oo n x))
= (Dv_k \oo 0 (x\oon n x)) = 0 .
\]

Hence,
\begin{equation}\label{eq:TestElement}
((v_k\oon 0 x)\oon n x) = 0 \in W
\end{equation}
for all $k\ge 0$.

On the other hand,
\[
((v_k\oon 0 x)\oon \lambda  x) = (\partial^kv_k\oon \lambda x)
= (-\lambda )^k (\partial+\lambda )^k v_k \ne 0,
\]
and if $k$ is large enough then the coefficient at  $\lambda ^n$ is nonzero,
a contradiction.
\end{proof}

\end{document}